\documentclass[a4paper,11pt]{article}
\usepackage[margin=1in]{geometry}
\usepackage{amsmath,amssymb,amsthm,mathtools}
\let\geq\geqslant
\let\leq\leqslant
\usepackage{lmodern,microtype}
\DeclarePairedDelimiter{\floor}{\lfloor}{\rfloor}
\DeclarePairedDelimiter{\ceil}{\lceil}{\rceil}
\usepackage[dvipsnames]{xcolor}
\usepackage[colorlinks=true,linkcolor=blue!50!black,citecolor=blue!50!black,urlcolor=blue!50!black]{hyperref}
\usepackage{enumitem}
\usepackage{sectsty}
\usepackage[longnamesfirst,numbers,sort&compress]{natbib}
\makeatletter
\def\NAT@spacechar{~}
\makeatother
\usepackage{thmtools,thm-restate}
\usepackage[noabbrev,capitalise]{cleveref}
\crefname{thm}{Theorem}{Theorems}
\Crefname{thm}{Theorem}{Theorems}
\crefname{lem}{Lemma}{Lemmas}
\Crefname{lem}{Lemma}{Lemmas}
\crefname{cor}{Corollary}{Corollaries}
\Crefname{cor}{Corollary}{Corollaries}
\crefname{prop}{Proposition}{Propositions}
\Crefname{prop}{Proposition}{Propositions}
\crefname{conj}{Conjecture}{Conjectures}
\Crefname{conj}{Conjecture}{Conjectures}
\crefname{question}{Question}{Questions}
\Crefname{question}{Question}{Questions}

\crefname{definition}{Definition}{Definitions}
\Crefname{definition}{Definition}{Definitions}
\crefname{rem}{Remark}{Remarks}
\Crefname{rem}{Remark}{Remarks}
\crefname{equation}{Equation}{Equations}
\Crefname{equation}{Equation}{Equations}

\makeatletter
\def\NAT@spacechar{~}
\makeatother
\allsectionsfont{\boldmath}
\newcommand{\mathdefn}[1]{\textcolor{Maroon}{#1}}
\newcommand{\defn}[1]{\textcolor{Maroon}{\emph{#1}}}
\newcommand{\had}{\operatorname{had}}
\newcommand{\hadodd}{\operatorname{had}_{\mathrm{o}}}
\newcommand{\haddom}{\operatorname{had}_{\mathrm{d}}}
\newcommand{\hadod}{\operatorname{had}_{\mathrm{od}}}
\newcommand{\hadsod}{\operatorname{had}_{\mathrm{sod}}}
\newcommand{\hadred}{\operatorname{had}_{\mathrm{rd}}}
\newcommand{\doi}[1]{\href{https://doi.org/#1}{\nolinkurl{doi:#1}}}

\newtheorem{thm}{Theorem}
\newtheorem{lem}[thm]{Lemma}
\newtheorem{cor}[thm]{Corollary}
\newtheorem{prop}[thm]{Proposition}
\newtheorem{conj}[thm]{Conjecture}

\theoremstyle{definition}

\title{\bfseries Independence number and\\ odd dominating models}
\author{Ken-ichi Kawarabayashi\footnote{National Institute of Informatics and the University of Tokyo, Japan (\texttt{k\_keniti@nii.ac.jp}). Research supported by Alexander von Humboldt Foundation (Humboldt Research Award), JSPS Kakenhi 26K21777 and JP25K24465, and by JST ASPIRE JPMJAP2302.}
\qquad David R. Wood\footnote{School of Mathematics, Monash University, Melbourne, Australia (\texttt{david.wood@monash.edu}). Research supported by the Australian Research Council.}
\qquad Jung Hon Yip\footnote{School of Mathematics, Monash University, Melbourne, Australia (\texttt{Junghon.Yip@monash.edu}). Supported by a Monash Graduate Scholarship.}}
\date{\today}
\usepackage[bottom]{footmisc}
\begin{document}
\maketitle
\begin{abstract}
    A \defn{dominating $K_t$-model} in a graph $G$ is a sequence $(T_1,\dots,T_t)$ of pairwise vertex-disjoint subtrees of $G$ such that for $1\leqslant i<j\leqslant t$ every vertex in $T_j$ has a neighbour in $T_i$. Replacing ``every vertex in $T_j$'' by ``some vertex in $T_j$'' recovers the standard definition of $K_t$-model, which is equivalent to a $K_t$-minor in~$G$. The \defn{dominating Hadwiger number} $\haddom(G)$ is the largest integer $t$ such that $G$ contains a dominating $K_t$-model. This notion was introduced by Illingworth and Wood [\emph{J. Graph Theory}, 2025].
    We prove two dominating-model analogues of improvements to the classical Duchet--Meyniel bound [\emph{Ann. Discrete Math.}, 1982].
    The first extends a result of Kawarabayashi, Plummer, and Toft [\emph{J. Combin. Theory Ser. B}, 2005].
    Every $n$-vertex graph $G$ with $\alpha(G)\geqslant2$ satisfies $\haddom(G)\geqslant(n+\omega(G))/(2\alpha(G)-1)$, where $\alpha(G)$ and $\omega(G)$ denote the independence number and clique number of a graph $G$ respectively.
    Second, we prove a constant-factor improvement of the Duchet--Meyniel bound for dominating models: every graph $G$ with $n\geq 1$ vertices satisfies $\haddom(G)  > 0.51075n/\alpha(G)$. This qualitatively extends the results of Fox [\emph{SIAM J. Discrete Math.}, 2010] and Balogh and Kostochka [\emph{Discrete Math.}, 2011].

    An \defn{odd dominating $K_t$-model} is a dominating $K_t$-model together with an assignment of one of two colours to each vertex of $V(T_1)\cup\dots\cup V(T_t)$, such that every edge of each $T_i$ is bichromatic and, for $1\leqslant i<j\leqslant t$, some edge between $T_i$ and $T_j$ is monochromatic. Let \defn{$\hadod(G)$} denote the largest integer $t$ such that $G$ has an odd dominating $K_t$-model. Our constant-factor improvement of the Duchet--Meyniel bound relies on our next result: every non-empty claw-free graph $G$ satisfies $\hadod(G) \geq 2\chi(G)/3$, generalising a result of Chudnovsky and Fradkin [\emph{J. Graph Theory}, 2010]. This $2/3$ constant is asymptotically best possible by a construction of K\"uhn, Sauermann, Steiner, and Wigderson [\emph{arXiv}:2512.20392].

    We also prove that, for $t\geqslant2$, every $n$-vertex graph $G$ with $n\geqslant1$ and no odd dominating $K_t$-model has an independent set of size at least $(n+t-1)/(2t-2)$. Equivalently, every $n$-vertex graph $G$ with $n\geqslant1$ satisfies $\hadod(G)\geqslant n/(2\alpha(G)-1)$.
    For the stronger notion in which every vertex of $T_j$ has a neighbour of the same colour in $T_i$, we prove the analogous results with $(n+t-1)/(2t-2)$ replaced by $(n+t-1)/(3t-3)$ and $2\alpha(G)-1$ replaced by $3\alpha(G)-1$.
\end{abstract}

\section{Introduction}\label{sec:intro}
Hadwiger's Conjecture~\citep{Hadwiger43} has long been considered one of the deepest problems in graph theory\footnote{All graphs in this paper are finite and simple. For a graph $G$, let $\alpha(G)$, $\omega(G)$, and $\Delta(G)$ denote its independence number, clique number, and maximum degree, respectively. A \defn{proper colouring} of $G$ assigns a colour to each vertex so that adjacent vertices receive distinct colours. The \defn{chromatic number} $\chi(G)$ is the minimum number of colours in a proper colouring of $G$. A graph $H$ is a \defn{minor} of a graph $G$ if a graph isomorphic to $H$ can be obtained from a subgraph of $G$ by contracting edges. A \defn{$K_t$-model} in $G$ is a sequence $(T_1,\dots,T_t)$ of pairwise vertex-disjoint subtrees of $G$, called \defn{branch sets}, such that for $1\leqslant i<j\leqslant t$ some vertex of $T_j$ has a neighbour in $T_i$. The complete graph $K_t$ is a minor of $G$ if and only if $G$ has a $K_t$-model. The \defn{Hadwiger number} $\had(G)$ is the maximum integer $t$ such that $G$ has a $K_t$-model. For any undefined graph-theoretic notation, see \citet{Diestel05}.}, since it generalises the well-known four-colour theorem \citep{AH77-I,AH77-II,RSST97} and suggests a deep relationship between graph minors and graph colouring; see \citep{Kawa15,Toft96,SeymourHC,CV20,Norin23,CLWY25} for surveys.
\begin{conj}[Hadwiger's Conjecture]
    \label{conj:hc}
    For every graph $G$, $\had(G) \geq \chi(G)$. Equivalently, if $G$ is $K_t$-minor-free, then $G$ is $(t-1)$-colourable.
\end{conj}

Hadwiger's Conjecture and several of its strengthenings have recently been disproved~\citep{OpenAIHC26,KSSW25,IS26}, while its linear relaxation has been proved~\citep{NS26,OpenAILinearHC26,Wood26} (see \cref{sec:conj}).
Here we study bounds on the independence number for dominating, odd dominating, and strongly odd dominating models, which we define later.
Every $n$-vertex graph $G$ with $n\geqslant1$ satisfies $\chi(G) \geq n/\alpha(G)$. So Hadwiger's Conjecture would imply:
\begin{conj}
    \label{conj:hc_alpha}
    For every $n$-vertex graph $G$ with $n\geqslant1$, $\had(G) \geq n/\alpha(G)$.
\end{conj}

\cref{conj:hc_alpha} has also been disproved~\citep{OpenAIHC26}. Nevertheless, \citet{DM82} proved using an elegant argument that every $n$-vertex graph $G$ with $n\geqslant1$ satisfies $\had(G) \geqslant n/(2\alpha(G)-1)$, which is within a factor of~2 of the bound predicted by \cref{conj:hc_alpha}. We shall refer to this as the \defn{Duchet--Meyniel bound}. For graphs with independence number~2, no constant-factor improvement of the resulting bound $\had(G)\geqslant n/3$ is known. More specifically, it is open if there is an absolute constant $\varepsilon>0$ such that every $n$-vertex graph $G$ with $\alpha(G)=2$ satisfies $\had(G)\geqslant(1/3+\varepsilon)n$. This question, due to \citet{SeymourHC}, has been the focus of much research~\citep{PST03,Fox10,BK11,CS12,SeymourHC}.
The recent counterexamples of \citet{OpenAIHC26} have independence number~2, and in this class the constant $1/3$ cannot be replaced by any constant greater than $26/75$.

Several improvements to the Duchet--Meyniel bound are known. \citet{KPT05} proved that every $n$-vertex graph $G$ with $\alpha(G)\geqslant2$ satisfies $\had(G)\geqslant(n+\omega(G))/(2\alpha(G)-1)$.
This bound was improved by \citet{PT10} to $\had(G)\geqslant(n+\omega(G)+2\alpha(G)-4)/(2\alpha(G)-1)$ for every $n$-vertex graph $G$ with $n\geqslant1$ that is neither complete nor edgeless.
For graphs with $\alpha(G)=2$, \citet{CS12} established the stronger bound $\had(G)\geqslant n/2$, provided that $\omega(G)\geqslant n/4$ when $n$ is even, or $\omega(G)\geqslant(n+3)/4$ when $n$ is odd.

For $\alpha(G)\geqslant3$, the bound $\had(G)\geqslant n/(2\alpha(G)-1.5)$ of \citet{KPT05} was improved to $\had(G)\geqslant n/(2\alpha(G)-2)$ by \citet{KS05}.
\citet{Wood-DMTCS07} obtained the related bound $\had(G)\geqslant(n+5\alpha(G)-5)/(2\alpha(G)-1)$ when $\had(G)\geqslant5$.
\citet{BLW11} then proved that every $n$-vertex graph $G$ with $n\geqslant1$ satisfies
\[
    \had(G)\geqslant
    \frac{n}{2\alpha(G)-\ceil{\log_{\tau}(\tau\alpha(G)/2)}},
    \qquad \tau:=\frac{2\sqrt{2}}{\sqrt{2}-1},
\]
which improves the bound of \citet{KS05} when $\alpha(G)\geqslant14$.

\citet{Fox10} gave the first improvement of the Duchet--Meyniel bound by an absolute constant factor: $\had(G)>0.504\,n/\alpha(G)$ for every $n$-vertex graph $G$ with $n\geqslant1$. This is asymptotically stronger than the bound of \citet{BLW11} when $\alpha(G)$ is large. \citet{BK11} then improved Fox's constant, proving that every $n$-vertex graph $G$ with $n\geqslant1$ satisfies $\had(G)\geqslant n/((2-c)\alpha(G))$ for a constant $c>1/19.2$.

We consider the following four strengthenings of $K_t$-models. In~(a), (c) and~(d), a \defn{colouring} is any function $c:\bigcup_iV(T_i)\to\{\text{red},\text{blue}\}$ (not necessarily a proper colouring of $G$). An edge $vw$ is \defn{bichromatic} if $c(v)\neq c(w)$, and \defn{monochromatic} if $c(v)=c(w)$.
\begin{enumerate}[label=(\alph*)]
    \item An \defn{odd $K_t$-model} (Gerards and Seymour, see~\citep{JT95}) is a $K_t$-model $(T_1,\dots,T_t)$ together with a colouring $c$ such that every edge of every $T_i$ is bichromatic and, for $1\leqslant i<j\leqslant t$, some edge of $G$ between $T_i$ and $T_j$ is monochromatic.
    \item A \defn{dominating $K_t$-model} (\citet{IW24}) is a sequence $(T_1,\dots,T_t)$ of pairwise vertex-disjoint subtrees of $G$ such that for $1\leqslant i<j\leqslant t$ every vertex of $T_j$ has a neighbour in $T_i$.
    \item An \defn{odd dominating $K_t$-model} is a dominating $K_t$-model $(T_1,\dots,T_t)$ together with a colouring $c$ such that every edge of every $T_i$ is bichromatic and, for $1\leqslant i<j\leqslant t$, some edge of $G$ between $T_i$ and $T_j$ is monochromatic. That is, $((T_1,\dots,T_t),c)$ is an odd $K_t$-model and $(T_1,\dots,T_t)$ is a dominating $K_t$-model.
    \item A \defn{strongly odd dominating $K_t$-model} is a dominating $K_t$-model $(T_1,\dots,T_t)$ together with a colouring $c$ such that every edge of every $T_i$ is bichromatic and, for $1\leqslant i<j\leqslant t$, every vertex $v$ of $T_j$ has a neighbour $u$ in $T_i$ with $c(u)=c(v)$. That is, the domination is witnessed by monochromatic edges.
\end{enumerate}
For each type of model, we define a Hadwiger-type parameter. Let $\hadodd(G)$, $\haddom(G)$, $\hadod(G)$ and $\hadsod(G)$ be the maximum integer $t$ such that $G$ has an odd, dominating, odd dominating and strongly odd dominating $K_t$-model, respectively. Every strongly odd dominating $K_t$-model is an odd dominating $K_t$-model, and every odd dominating $K_t$-model is an odd $K_t$-model and a dominating $K_t$-model. Thus every graph $G$ satisfies
\[
    \omega(G)\leqslant\hadsod(G)\leqslant\hadod(G)
    \leqslant\hadodd(G)\leqslant\had(G), \, \text{and}
\]
\[
    \omega(G)\leqslant\hadsod(G)\leqslant\hadod(G)
    \leqslant\haddom(G)\leqslant\had(G).
\]

In some ways, dominating models resemble ordinary complete graph models. For $t\leqslant4$, \citet{IW24} proved that every graph with no dominating $K_t$-model is $(t-1)$-colourable. The case $t=5$ was proved by \citet{D4CT26}: every graph with no dominating $K_5$-model is $4$-colourable, thus extending the four-colour theorem. Moreover, \citet{IW24} proved that $\haddom(G)\geqslant\chi(G)$ for almost every graph\footnote{Here, `almost every' means that for each $p \in (0, 1)$, a random graph $G \sim \mathcal{G}(n, p)$ satisfies the statement with probability tending to one as $n\to\infty$.}, mirroring the corresponding result of \citet{BCE80} for ordinary models. Graphs with no dominating $K_t$-model are also reasonably sparse: \citet{GNTT26} proved that graphs with no dominating $K_t$-model are $O(t(\log t)^2)$-degenerate, while classical results of \citet{Kostochka82,Kostochka84} and \citet{Thomason84} show that $K_t$-minor-free graphs are $O(t\sqrt{\log t})$-degenerate.

Our first goal is to extend two refinements of the Duchet--Meyniel bound to dominating models, adding to the growing collection of parallel results between ordinary and dominating models. As mentioned above,
\citet[Theorem~1.2]{KPT05} proved that every graph $G$ with $\alpha(G)\geqslant2$ satisfies $\had(G)\geqslant (|V(G)|+\omega(G))/(2\alpha(G)-1)$, strengthening the Duchet--Meyniel bound by incorporating the clique number. In the dominating setting, \citet{IW24} proved that $\haddom(G)\geqslant |V(G)|/(2\alpha(G)-1)$ for every graph $G$. We extend both these results.
\begin{restatable}{thm}{cliquenumberthm}\label{thm:omega}
    Every $n$-vertex graph $G$ with $\alpha(G)\geqslant2$ satisfies
    \[ \haddom(G)\ \geqslant\ \frac{n+\omega(G)}{2\alpha(G)-1}. \]
\end{restatable}
Our next result gives a constant-factor improvement of the Duchet--Meyniel bound in the dominating setting, which qualitatively extends the results of \citet{Fox10} and \citet{BK11}.
\begin{restatable}{thm}{uniformimprovementthm}
    \label{thm:uniform_improvement}
    Let $c:= (49 - \sqrt{2113})/72$. Every $n$-vertex graph $G$ with $n\geqslant1$ satisfies
    \begin{equation*}
        \haddom(G) \geq \frac{n}{(2-c)\alpha(G)}
        = \frac{72}{95+\sqrt{2113}}\frac{n}{\alpha(G)}
        > 0.51075\frac{n}{\alpha(G)}.
    \end{equation*}
\end{restatable}
Our proof follows the ideas of \citet{Fox10} and \citet{BK11}, but uses a simpler analysis at the expense of a weaker constant.

A key additional ingredient in the proof of \cref{thm:uniform_improvement} is the following result on claw-free graphs.
A \defn{claw} is an induced $K_{1,3}$, and a graph is \defn{claw-free} if no induced subgraph is isomorphic to a claw. \citet{CF10} proved that every non-empty claw-free graph $G$ satisfies $\had(G) \geq 2\chi(G)/3$, proving an approximate version of Hadwiger's Conjecture for claw-free graphs. We strengthen this to odd dominating models.
\begin{restatable}{thm}{clawfreethmodddom}
    \label{thm:odd_dom-claw-free}
    Every non-empty claw-free graph $G$ satisfies $\hadod(G) \geq 2\chi(G) / 3$.
\end{restatable}

The proof of the original claw-free result of \citet{CF10} relies on a detailed analysis using the structure theorem for claw-free graphs~\citep{CS05}. Subsequent results \citep{CS12,CS10,KR15} allow us to give a shorter proof that generalises the bound to odd dominating models. Our proof of \cref{thm:odd_dom-claw-free} still relies indirectly on the structure theorem of claw-free graphs since we use a bound of \citet{KR15} (see \cref{thm:reed}).

The Duchet--Meyniel theorem has been extended to settings~(a) and~(b). In setting~(a), \citet{KS07} proved that every $n$-vertex graph $G$ with $n\geqslant1$ satisfies $\hadodd(G)\geqslant n/(2\alpha(G)-1)$. In setting~(b), \citet[Theorem~14]{IW24} proved that for every graph $G$ with $n\geqslant1$ vertices there is an integer $h$ with $1\leqslant h\leqslant\haddom(G)$ such that $G$ has an $h$-colourable induced subgraph on at least $(n+h)/2$ vertices. Hence $\haddom(G)\geqslant n/(2\alpha(G)-1)$.
Our next result extends the Duchet--Meyniel theorem to setting~(c).

\begin{thm}\label{thm:A}
    For every $n$-vertex graph $G$ with $n\geqslant1$, there is an integer $h$ with $1\leqslant h\leqslant\hadod(G)$ such that $G$ has an $h$-colourable induced subgraph on at least $(n+h)/2$ vertices. In particular,
    \[ \hadod(G)\ \geqslant\ \frac{n}{2\alpha(G)-1}. \]
\end{thm}

Since $\hadod(G)\leqslant\min\{\hadodd(G),\haddom(G)\}$, \cref{thm:A} implies the above results of \citet{DM82}, \citet{KS07}, and \citet{IW24}. In \cref{sec:A} we prove a slightly stronger result (\cref{thm:Aplus}).

Does the Duchet--Meyniel bound also hold for strongly odd dominating models? The answer is no\footnote{
    Consider an $n$-vertex graph $G$ with independence number $2$ with clique number $\omega(G) = \Theta(\sqrt{n \log n})$; such graphs exist \citep{Kim95}. Let $(T_1, \dots, T_t)$ be a strongly odd dominating $K_t$-model with a colouring $c$. In any colouring, all singleton branch sets have the same colour and each non-singleton branch set contains both colours. By swapping red and blue, we may ensure that blue occurs in every branch set. For $i \geq 1$, let $\mathcal{B}_i$ be the collection of branch sets of size $i$. Let $\mathcal{S}_r$ and $\mathcal{S}_b$ denote the collections of branch sets with size $3$ whose degree-$2$ vertex is coloured red and blue, respectively. Then $\mathcal{B}_3:= \mathcal{S}_r \cup \mathcal{S}_b$. Consider the strongly odd dominating model induced by $\mathcal{S}_b \cup \mathcal{B}_2 \cup \mathcal{B}_1$. Since this model is strongly odd and dominating, the blue vertices form a clique. However, since $\omega(G) = \Theta(\sqrt{n \log n})$, we deduce that $|\mathcal{S}_b| +|\mathcal{B}_2| + |\mathcal{B}_1| = O(\sqrt{n \log n})$. Now consider the strongly odd dominating model induced by $\mathcal{S}_r$. Again, the red vertices form a clique, so $|\mathcal{S}_r| = O(\sqrt{n \log n})$.
    We deduce that $|\mathcal{B}_1 \cup \mathcal{B}_2 \cup \mathcal{B}_3| = O(\sqrt{n \log n})$. Let $s:=|\mathcal{B}_1 \cup \mathcal{B}_2 \cup \mathcal{B}_3|$.

    Since the branch sets are pairwise vertex-disjoint, there are $t-s$ branch sets containing at least four vertices. Therefore $n\geq s+4(t-s)=4t-3s$, and so $t \leq n/4 + 3s/4 = n/4 + O(\sqrt{n\log n})$, which is much less than the predicted bound of $n/3$.

}. Nevertheless, we show an analogue of \cref{thm:A}, with $(n+h)/2$ replaced by $(n+h)/3$ and $2\alpha(G)-1$ replaced by $3\alpha(G)-1$.

\begin{thm}\label{thm:B}
    For every graph $G$ with $n\geqslant1$ vertices, there is an integer $h$ with $1\leqslant h\leqslant\hadsod(G)$ such that $G$ has an $h$-colourable induced subgraph on at least $(n+h)/3$ vertices. In particular,
    \[ \hadsod(G)\ \geqslant\ \frac{n}{3\alpha(G)-1}. \]
\end{thm}

\textbf{Outline.}
We discuss related colouring problems in \cref{sec:conj}. We then prove \cref{thm:omega} in \cref{sec:omega}, followed by \cref{thm:odd_dom-claw-free} in \cref{s:claw-free}. Next we prove \cref{thm:uniform_improvement}
in \cref{s:uniform_improvement}. Finally, \cref{sec:A,sec:B} establish the odd dominating and strongly odd dominating bounds (\cref{thm:A,thm:B}), together with the corresponding upper bounds on chromatic number (\cref{prop:wagner,prop:sod}).

\subsection{Graph Minors and Colouring}\label{sec:conj}

We now place \cref{thm:omega,thm:uniform_improvement,thm:odd_dom-claw-free,thm:A,thm:B} in the context of more general conjectures and related problems studied in the literature.

Hadwiger's Conjecture holds for $t\leqslant4$~\citep{Hadwiger43,Dirac52}. \citet{Wagner37} proved that the $t=5$ case is equivalent to the four-colour theorem, and the $t=6$ case was proved by \citet*{RST-Comb93}, also using the four-colour theorem.

For every $K_t$-minor-free graph $G$, \citet{Kostochka84} and \citet{Thomason84} independently proved that $\chi(G)=O(t\sqrt{\log t})$. This was improved to $\chi(G)=O(t(\log t)^{1/4+\varepsilon})$ for every fixed $\varepsilon>0$ by \citet{NPS23}.
The bound was subsequently improved to $\chi(G)=O(t\log\log t)$ by \citet{DP24}, to $\chi(G)=O(t\sqrt{\log\log t})$ by \citet{Lin26}, and to $\chi(G)=O(t\log\log\log t)$ by \citet{LL26}.
Recently, \citet{NS26} proved the linear Hadwiger Conjecture: every $K_t$-minor-free graph satisfies $\chi(G)\leq Ct$ for an absolute constant $C$. \citet{Wood26} obtained an explicit bound $\chi(G)\leq 1{,}118{,}430\,t$. OpenAI also announced a qualitatively stronger linear list-colouring bound: $\chi_{\mathrm{list}}(G)\leq Ct$ for every $K_t$-minor-free graph $G$~\citep{OpenAILinearHC26}, where $\chi_{\mathrm{list}}(G)$ is the least $k$ such that $G$ can be properly coloured from any assignment of a list of $k$ available colours to each vertex.

The case of Hadwiger's Conjecture for graphs with independence number $2$ has attracted particular attention~\citep{PST03,CS12,Blasiak07,CO08,Bosse19,NS26a,CLWY25}. \citet{PST03} proved that, for graphs with independence number~2, Hadwiger's Conjecture is equivalent to \cref{conj:hc_alpha}. \citet{SeymourHC} writes ``My own belief is, if Hadwiger's Conjecture is true for graphs with stability number two then it is probably true in general, so it would be very nice to decide this case''. Every graph of independence number $2$ is claw-free.
However, Hadwiger's Conjecture is false, shown recently by \citet{OpenAIHC26}, and their counterexamples have independence number $2$. They showed:
\begin{thm}[{\citet[Corollary 1.2]{OpenAIHC26}}]
    \label{thm:disproof_hc}
    There are $n$-vertex graphs $G$ of independence number $2$, for arbitrarily large $n$, for which
    \[
        \chi(G)>\frac{75}{52}\had(G)-\frac{25}{26}.
    \]
    Thus Hadwiger's Conjecture is false.
\end{thm}

Another strengthening of Hadwiger's Conjecture is the Odd Hadwiger Conjecture, proposed by Gerards and Seymour~\citep{SeymourHC}. It asserts that for every graph $G$, $\hadodd(G)\geq\chi(G)$.

The Odd Hadwiger Conjecture holds for $t\leq 4$, and \citet{Guenin05} announced a proof for $t=5$; see \citet{SeymourHC} for background. \citet*{GRRSV09} proved that every graph with no odd $K_t$-minor is $O(t\sqrt{\log t})$-colourable; \citet{Kawa09} subsequently gave a shorter proof. \citet{Steiner22} proved that if every $K_t$-minor-free graph is $f(t)$-colourable, then every graph with no odd $K_t$-minor is $2f(t)$-colourable. Together with the linear bound of \citet{NS26}, this implies that every graph with no odd $K_t$-minor satisfies $\chi(G)\leq Ct$ for some absolute constant $C$.

Nevertheless, \citet*{KSSW25} disproved the Odd Hadwiger Conjecture. Their counterexamples also have independence number $2$. Since these graphs are claw-free, their construction also shows that the constant $2/3$ in \cref{thm:odd_dom-claw-free} is asymptotically best possible.
\begin{thm}[{\citet[Theorem 1.3]{KSSW25}}]
    \label{thm:disproof_odd_hc}
    For every $\varepsilon>0$ and all sufficiently large $t$, there is a graph $G$ with independence number $2$ with no odd $K_t$-minor, and $\chi(G)\geq (3/2-\varepsilon)t$.
\end{thm}

\citet{IW24} asked whether $\chi(G)\leqslant\haddom(G)$ for every graph $G$. This question became known as the Dominating Hadwiger Conjecture. For $t\leqslant3$, a graph has a dominating $K_t$-model if and only if it has a $K_t$-model~\citep{IW24}. So $\haddom(G)=\had(G)$ whenever $\haddom(G)\leqslant2$, and the Dominating Hadwiger Conjecture holds for graphs $G$ with $\haddom(G)\leqslant2$. \citet{IW24} proved that every graph $G$ with $\haddom(G)\leqslant3$ is 3-colourable, and that $\chi(G)\leqslant 2^{\haddom(G)-1}$ for every graph $G$. \citet{D4CT26} proved that every graph $G$ with $\haddom(G)\leqslant4$ is 4-colourable. \citet{GNTT26} proved that every graph with average degree at least $Ct(\log t)^2$ has a dominating $K_t$-model, where $C$ is an absolute constant. Since a dominating $K_t$-model in a subgraph of $G$ is also a dominating $K_t$-model in $G$, it follows that every graph with no dominating $K_t$-model is $O(t(\log t)^2)$-degenerate, and is thus $O(t(\log t)^2)$-colourable. \citet{SS25} proved the Dominating Hadwiger Conjecture for several classes of graphs $G$ with $\alpha(G)\leqslant2$, and \citet{ST25} proved it for $2K_2$-free graphs. \citet{IS26} disproved the Dominating Hadwiger Conjecture, using graphs $G$ with $\alpha(G)\leqslant2$. Since $\haddom(G)\leqslant\had(G)$, the OpenAI counterexample~\citep{OpenAIHC26} is also a counterexample to the Dominating Hadwiger Conjecture.

For odd dominating and strongly odd dominating models, the case $t=3$ is easy. For every graph $G$,\footnote{First, every bipartite graph $G$ satisfies $\hadodd(G)\leqslant2$. To see this, let $b$ be a proper 2-colouring of $G$, and let $((T_1,T_2,T_3),c)$ be an odd $K_3$-model in $G$. Take both $b$ and $c$ to use red and blue. Two proper 2-colourings of a tree are equal up to swapping the colours, so on each $T_i$ the colouring $c$ either agrees with $b$ or is obtained from $b$ by swapping red and blue. A monochromatic edge $uv$ with $u\in V(T_i)$ and $v\in V(T_j)$ satisfies $b(u)\neq b(v)$, so $c$ agrees with $b$ on exactly one of $T_i$ and $T_j$. This is impossible for all three pairs of branch sets. Second, every non-bipartite graph $G$ has a strongly odd dominating $K_3$-model. To see this, let $v_0v_1\dots v_{2k}v_0$ be an odd cycle in $G$. Let $T_1$ be the path $v_0v_1\dots v_{2k-2}$, coloured alternately with $v_0$ coloured red (so $v_{2k-2}$ is also coloured red), let $T_2:=\{v_{2k-1}\}$ and $T_3:=\{v_{2k}\}$, and colour $v_{2k-1}$ and $v_{2k}$ with colour red. Then $v_{2k-1}$ is adjacent to $v_{2k-2}$, and $v_{2k}$ is adjacent to $v_0$ and $v_{2k-1}$, all of which have colour red, so $((T_1,T_2,T_3),c)$ is a strongly odd dominating $K_3$-model.}
\[ \hadsod(G)\leqslant2 \iff \hadod(G)\leqslant2 \iff G\text{ is bipartite}
    . \]

Consequently, $\chi(G)\leqslant\hadod(G)$ for every graph $G$ with $\hadod(G)\leqslant2$, and $\chi(G)\leqslant\hadsod(G)$ for every graph $G$ with $\hadsod(G)\leqslant2$. Further, for every non-empty graph $G$, \cref{prop:wagner,prop:sod} give the exponential bounds $\chi(G)\leqslant2^{\hadod(G)-1}$ and $\chi(G)\leqslant3\cdot2^{\hadsod(G)}-5$.

In light of \citet{OpenAIHC26}'s construction showing that $\had(G)/\chi(G)$ can be asymptotically at most $52/75 = 2/3 + 2/75$, together with the recent proof of the linear Hadwiger Conjecture \citep{NS26,OpenAILinearHC26,Wood26}, it is natural to ask for the best possible ratio between the Hadwiger number and chromatic number for a graph class. For a class $\mathcal{G}$ of graphs, define
\begin{equation*}
    \rho(\mathcal{G})
    :=\inf\left\{\frac{\had(G)}{\chi(G)}:
    G\in\mathcal{G} \text{ is non-empty}\right\},
\end{equation*}
and define $\rho_d(\mathcal{G})$ and $\rho_o(\mathcal{G})$ analogously using $\haddom(G)$ and $\hadodd(G)$, respectively.
For the class of claw-free graphs, we write $\rho:=\rho(\text{claw-free})$, $\rho_d:=\rho_d(\text{claw-free})$, and $\rho_o:=\rho_o(\text{claw-free})$.
For claw-free graphs, \citet{CF10} and \cref{thm:odd_dom-claw-free,thm:disproof_hc} give
\[
    \frac23\leq \rho_d\leq \rho \leq \frac23+\frac{2}{75},
\]
where the upper bound follows from \citet{OpenAIHC26}.
Moreover, \cref{thm:odd_dom-claw-free,thm:disproof_odd_hc} give $\rho_o=2/3$, with the upper bound due to \citet{KSSW25}.

\section{Dominating models and the clique number}\label{sec:omega}
We need the following preliminary results.
\begin{lem}[\citet{IW24}]
    Every $n$-vertex graph $G$ with $n \geq 1$ satisfies
    \begin{equation}\label{eq:DM}
        n\leqslant(2\alpha(G)-1)\,\haddom(G).
    \end{equation}
\end{lem}

We also use the following lemma of \citet[Lemma~2.5]{FP21}; see also \citet[10.1]{CLSS17} for a short proof.

\begin{lem}[{\citet{FP21}}]
    \label{lem:FP}
    Every non-empty connected graph $H$ has a set $S\subseteq V(H)$ such that $H[S]$ is connected, $|S|\leqslant 2\alpha(H)-1$ and $\alpha(H[S])=\alpha(H)$.
\end{lem}

We now prove \cref{thm:omega}, following the proof of Theorem~1.2 of \citet{KPT05}.
Using \cref{lem:FP} yields a slightly shorter argument than the original proof of \citet{KPT05}.

\begin{proof}[Proof of \cref{thm:omega}]
    Equivalently, we show that every $n$-vertex graph $G$ with $\alpha(G)\geqslant2$ satisfies
    \[ n+\omega(G)\ \leqslant\ (2\alpha(G)-1)\,\haddom(G). \]
    We proceed by induction on $n$. Let $a:=\alpha(G)\geqslant2$ and $h:=\haddom(G)$. Let $Q$ be a maximum clique of $G$, and let $H:=G-Q$. Then $|Q|=\omega(G)\leqslant h$.

    \emph{Case 1: $\alpha(H)\leqslant a-1$.} Since $a\geqslant2$, the graph $G$ is not complete, so $H$ is non-empty, and \cref{eq:DM} applied to $H$ gives $|V(H)|\leqslant(2\alpha(H)-1)\haddom(H)\leqslant(2a-3)h$. Hence
    \[ n+\omega(G)=|V(H)|+2|Q|\leqslant(2a-3)h+2h=(2a-1)h. \]

    \emph{Case 2: $\alpha(H)=a$ and $H$ is connected.} By \cref{lem:FP}, there is a set $S\subseteq V(H)$ such that $G[S]$ is connected, $|S|\leqslant2a-1$ and $\alpha(G[S])=a$. Every vertex $v\in V(G)\setminus S$ has a neighbour in $S$, since otherwise $v$ together with a maximum independent set of $G[S]$ is an independent set of size $a+1$ in $G$. Let $G':=G-S$. If $(T_1,\dots,T_k)$ is a dominating $K_k$-model in $G'$, then $(T_0,T_1,\dots,T_k)$ is a dominating $K_{k+1}$-model in $G$, where $T_0$ is a spanning tree of $G[S]$. Thus $\haddom(G')\leqslant h-1$. Since $Q\subseteq V(G')$, we have $\omega(G')=|Q|$, and $G'$ is non-empty. If $\alpha(G')\geqslant2$, then by induction, $|V(G')|+|Q|\leqslant(2\alpha(G')-1)\haddom(G')$. If $\alpha(G')=1$, then $G'$ is complete. Since $Q$ is a maximum clique of $G$ and $Q \subseteq V(G')$, we have $V(G') = Q$, so $|V(G')|+|Q|=2|V(G')|=2\haddom(G')$. Since $\alpha(G')\leqslant a$ and $a\geqslant2$, in both cases
    \[ |V(G')|+|Q|\leqslant(2a-1)(h-1). \]
    Adding $|S|\leqslant 2a-1$ gives $n+\omega(G)\leqslant(2a-1)h$.

    \emph{Case 3: $\alpha(H)=a$ and $H$ is disconnected.} Partition $V(H)$ into non-empty sets $A_1$ and $A_2$ with no edge between $A_1$ and $A_2$, and let $a_i:=\alpha(H[A_i])$ for $i\in\{1,2\}$. Then $a_1+a_2=a$ and $a_1,a_2\geqslant1$. For each $q\in Q$ there exists $i\in\{1,2\}$ such that $q$ has a neighbour in every maximum independent set of $H[A_i]$. Otherwise, for each $i\in\{1,2\}$ there is a maximum independent set $J_i$ of $H[A_i]$ containing no neighbour of $q$, and then $J_1\cup J_2\cup\{q\}$ is an independent set of size $a+1$ in $G$, which is a contradiction. Assign each $q\in Q$ to one such $i$, let $Q_i$ be the set of vertices of $Q$ assigned to $i$, and let $G_i:=G[A_i\cup Q_i]$.

    We claim that $\alpha(G_i)=a_i$. Since $H[A_i]$ is an induced subgraph of $G_i$, we have $\alpha(G_i)\geqslant a_i$. For the reverse inequality, let $X$ be an independent set in $G_i$. Since $Q_i$ is a clique, $|X\cap Q_i|\leqslant1$. If $X\cap Q_i=\varnothing$, then $|X|\leqslant a_i$. If $X\cap Q_i=\{q\}$, then $X\setminus\{q\}$ is an independent set of $H[A_i]$ containing no neighbour of $q$, so $X\setminus\{q\}$ is not a maximum independent set of $H[A_i]$, and $|X|\leqslant(a_i-1)+1=a_i$. This proves the claim.

    Moreover, $|Q_i|\leqslant\omega(G_i)$ and $\haddom(G_i)\leqslant h$, and
    \[ n+|Q|=\sum_{i=1}^2\big(|V(G_i)|+|Q_i|\big). \]
    If $a=2$, then $a_1=a_2=1$, so each $G_i$ is complete and $|V(G_i)|\leqslant|Q|$. Hence $n+|Q|=|V(G_1)|+|V(G_2)|+|Q|\leqslant3|Q|\leqslant3h=(2a-1)h$.

    Now assume that $a\geqslant3$. Each $G_i$ has fewer vertices than $G$ (since $A_{3-i}\neq\varnothing$), and
    \[ n+|Q|\leqslant\sum_{i=1}^2\big(|V(G_i)|+\omega(G_i)\big). \]
    If $a_1,a_2\geqslant2$, then by induction the right-hand side is at most $(2a_1-1)h+(2a_2-1)h=(2a-2)h$. Otherwise, without loss of generality, $a_1=1$ and $a_2=a-1\geqslant2$. Then $G_1$ is complete, so $|V(G_1)|+\omega(G_1)=2|V(G_1)|=2\haddom(G_1)\leqslant 2h$, and $|V(G_2)|+\omega(G_2)\leqslant(2a-3)h$ by induction. Thus the right-hand side is at most $(2a-1)h$. In both cases, $n+\omega(G)\leqslant(2a-1)h$.
\end{proof}

\section[Proof of the claw-free theorem]{Claw-free graphs}
\label{s:claw-free}
We need the following results. The first result is a straightforward extension of a deep result of \citet[1.3]{CS12} to odd dominating models. For completeness, we provide the proof in the \hyperref[sec:appendix]{Appendix}.
We remark that related results for $\haddom$ and $\hadodd$ were previously proved by \citet[Theorem 1.6]{SS25} and \citet[Theorem 1.7]{JSWZ25}, respectively.

\begin{restatable}{lem}{dominatingseagulls}
    \label{thm:dominating_seagulls}
    Every $n$-vertex graph $G$ with $n\geqslant1$, $\alpha(G) \leq 2$ and
    $2\omega(G) \geq \ceil{n/2} + 1$ has
    $\hadod(G) \geq \ceil{n/2}$.
\end{restatable}

A graph $G$ is \defn{tame} if there is a connected claw-free graph $H$ with $\alpha(H) \geq 3$ such that $G$ is an induced subgraph of $H$.
\begin{thm}[{\citet[1.7]{CS10}}]
    \label{thm:tame}
    Let $G$ be tame. Then either $\alpha(G) \geq 3$, or $V(G)$ is the union of four cliques.
\end{thm}
In particular, if $G$ is tame and $\alpha(G)\leq 2$, then \cref{thm:tame} implies that $V(G)$ is the union of four cliques, so $\omega(G)\geq |V(G)|/4$. This linear bound is stronger than what Ramsey theory gives: there are graphs with $\alpha(G)=2$ and $\omega(G)=\Theta(\sqrt{|V(G)|\log |V(G)|})$ \citep{Kim95}.

The next result concerns Reed's $(\omega, \Delta, \chi)$-conjecture. \citet{Reed98} conjectured that every non-empty graph $G$ satisfies $\chi(G) \leq \ceil{(\Delta(G) + 1 + \omega(G))/2}$. The conjecture remains open in general, but it was proved for claw-free graphs by \citet{KR15} using the structure theorem of claw-free graphs~\citep{CS05}.

\begin{thm}[{\citet[Theorem 1.6]{KR15}}]
    \label{thm:reed}
    If $G$ is a non-empty claw-free graph, then $\chi(G) \leq \ceil{(\Delta(G) + 1 + \omega(G))/2}$.
\end{thm}

A graph $G$ is \defn{vertex-critical} if $\chi(G-v)<\chi(G)$ for every $v\in V(G)$. A graph $G$ is \defn{decomposable} if $V(G)$ can be partitioned into two subsets $A$ and $B$ with $A,B\neq\varnothing$ such that each vertex of $A$ is adjacent to each vertex of $B$. Equivalently, $G$ is decomposable if the complement $\overline{G}$ is disconnected. We need the following deep theorem of \citet{gallai_kristische_1963} (see \citet{stehlik_critical_2003} for a strengthening).

\begin{thm}[\citet{gallai_kristische_1963}]
    \label{thm:deep_theorem_of_gallai}
    Every vertex-critical graph $G$ with $|V(G)| < 2\chi(G) - 1$ is decomposable.
\end{thm}

Lastly, we need the following modification of a result of \citet*[Theorem 3.3]{PST03} to odd dominating models. The proof of \cref{thm:plummer} is a straightforward modification of the proof in \citet[Theorem 3.3]{PST03}; we provide the proof in the \hyperref[sec:appendix]{Appendix}. Observe that \cref{thm:plummer} is a strengthening of \cref{thm:omega} in the $\alpha(G) \leq 2$ case.
\begin{restatable}[{\citet[Theorem 3.3(a)]{PST03}, modified}]{lem}{plummer}
    \label{thm:plummer}
    Every $n$-vertex graph $G$ with $n\geqslant1$
    and $\alpha(G) \leq 2$ satisfies
    $\hadod(G) \geq \ceil{(\omega(G) + n)/3}$.
\end{restatable}

We now prove \cref{thm:odd_dom-claw-free}.

\clawfreethmodddom*

\begin{proof}
    Let $G$ be a minimal claw-free counterexample, and let $n:=|V(G)|$. The graph $G$ is connected, since otherwise a component of maximum chromatic number is a smaller counterexample. Then $\hadod(G) < 2\chi(G) / 3$.
    If $\alpha(G) = 1$, then the result is immediate. First suppose that $\alpha(G) = 2$.
    Since $G$ is a minimal counterexample, for each vertex $v \in V(G)$, $\hadod(G - v) \geq 2\chi(G -v)/3$. Since $\hadod(G - v) \leq \hadod(G)$, we have $2\chi(G - v)/3 \leq \hadod(G - v) \leq \hadod(G) < 2\chi(G) / 3$. Thus $\chi(G - v) < \chi(G)$, and we deduce that $G$ is vertex-critical. If $G$ is decomposable, then $V(G) = A \cup B$, where $A \cap B = \varnothing$, $A \neq \varnothing$ and $B \neq \varnothing$, and each vertex of $A$ is adjacent to each vertex of $B$. Since $G$ is a minimal counterexample, $\hadod(G[A]) \geq 2\chi(G[A])/3$ and $\hadod(G[B]) \geq 2\chi(G[B])/3$. By swapping red and blue if necessary\footnote{In any odd clique model coloured red and blue, all singleton branch sets have the same colour and each non-singleton branch set contains both colours. By swapping red and blue, we may ensure that blue occurs in every branch set.}, choose the colourings of odd dominating models in $G[A]$ and $G[B]$ so that blue occurs in every branch set. Since every vertex of $A$ is adjacent to every vertex of $B$, these models can be concatenated, so $\hadod(G) \geq \hadod(G[A]) + \hadod(G[B])$. Since $\chi(G) = \chi(G[A]) + \chi(G[B])$, we obtain $\hadod(G) \geq 2\chi(G)/3$, a contradiction. Consequently $G$ is not decomposable, and \cref{thm:deep_theorem_of_gallai} implies that $n \geq 2\chi(G) - 1$.\footnote{In fact, $n=2\chi(G)-1$: since $\alpha(G)=2$, we have $n\leq 2\chi(G)$, while $n=2\chi(G)$ would give $\chi(G-v)\geq(n-1)/2>\chi(G)-1$ for every $v\in V(G)$, contradicting vertex-criticality since $\chi(G-v)$ is an integer. This argument is due to \citet{PST03}.} By \cref{thm:plummer}, $\hadod(G) \geq (n+1)/3$. Thus $3\hadod(G) -1 \geq n \geq 2\chi(G) - 1$ and $3\hadod(G) \geq 2\chi(G)$, as desired.

    Hence we may assume that $\alpha(G)\geq 3$.
    Let $v$ be a vertex of maximum degree and let $Z:=G[N_G(v)]$. Then $|V(Z)|=\Delta(G)$ and $\alpha(Z)\leq 2$, since $G$ is claw-free. Moreover, $Z$ is tame, since $G$ is connected and $\alpha(G)\geq 3$. By \cref{thm:tame}, $Z$ is the union of four cliques, so $\omega(Z)\geq\Delta(G)/4$.

    First suppose that $\hadod(G)$ is even. Choose a (possibly empty) clique $Q$ of $Z$ of size exactly $\floor{\Delta(G)/4}$, and let $q:= |Q|$. Let $Z'$ be an induced subgraph of $Z$ containing $Q$ such that $|V(Z')| = 4q$; this is possible since $4q \leq \Delta(G)$. Let $Z'':=G[V(Z')\cup\{v\}]$. Then $|V(Z'')|=4q+1$, $\alpha(Z'')\leq 2$, and $\omega(Z'')\geq q+1$, and consequently,
    \begin{equation*}
        2\omega(Z'')\geq 2q+2
        =\ceil{|V(Z'')|/2}+1.
    \end{equation*}
    By \cref{thm:dominating_seagulls} applied to $Z''$, $\hadod(G)\geq \hadod(Z'')\geq 2q+1$. Since $\hadod(G)$ is even, $\hadod(G)\geq 2q+2$. Therefore, $\Delta(G)\leq 4q+3\leq 2\hadod(G)-1$.

    Now suppose that $\hadod(G)$ is odd. Choose a clique $Q$ of $Z$ of size exactly $\ceil{\Delta(G)/4}$, and let $q:=|Q|$. Let $Z'$ be an induced subgraph of $Z$ containing $Q$ such that $|V(Z')|=4q-3$; this is possible since $4q-3\leq\Delta(G)$. Then $\alpha(Z')\leq 2$ and $\omega(Z')\geq q$, and consequently,
    \begin{equation*}
        2\omega(Z')\geq 2q
        =\ceil{|V(Z')|/2}+1.
    \end{equation*}
    By \cref{thm:dominating_seagulls} applied to $Z'$, $\hadod(Z')\geq 2q-1$. By swapping red and blue if necessary, we may choose a two-colouring of this model so that blue occurs in every branch set, and colour $v$ blue. Since every vertex of $Z'$ is adjacent to $v$, appending $\{v\}$ as the last branch set gives $\hadod(G)\geq \hadod(Z')+1\geq 2q$. Since $\hadod(G)$ is odd, $\hadod(G)\geq 2q+1$. Therefore, $\Delta(G)\leq 4q\leq 2\hadod(G)-2$.

    We deduce that
    \begin{equation}
        \label{eq:claw-free-degree-bounds}
        \Delta(G)\leq
        \begin{cases}
            2\hadod(G)-1, & \text{if }\hadod(G)\text{ is even}, \\
            2\hadod(G)-2, & \text{if }\hadod(G)\text{ is odd}.
        \end{cases}
    \end{equation}

    Finally, \cref{thm:reed,eq:claw-free-degree-bounds} and $\omega(G)\leq \hadod(G)$ give
    \[
        2\chi(G)\leq\Delta(G)+\omega(G)+2
        \leq
        \begin{cases}
            3\hadod(G)+1, & \text{if }\hadod(G)\text{ is even}, \\
            3\hadod(G),   & \text{if }\hadod(G)\text{ is odd}.
        \end{cases}
    \]
    Since $2\chi(G)$ is even, in either case $2\chi(G)\leq 3\hadod(G)$, contradicting the choice of $G$.
\end{proof}

Note that the proof of \cref{thm:odd_dom-claw-free} shows a stronger result: every non-empty claw-free graph $G$ has an odd dominating $K_t$-model with $t\geq 2\chi(G)/3$ in which every branch set has at most three vertices. This is because the models constructed in \cref{thm:plummer,thm:dominating_seagulls} have branch sets that are either induced three-vertex paths or single vertices.

One might wonder if \cref{thm:odd_dom-claw-free} can be strengthened further by replacing $\hadod$ with $\hadsod$. More precisely, does every claw-free graph $G$ satisfy $\hadsod(G)\geq 2\chi(G)/3$?
This is false\footnote{As mentioned, an $n$-vertex graph $G$ with independence number $2$ with clique number $\omega(G) = \Theta(\sqrt{n \log n})$ has $\hadsod(G) \leq n/4 + O(\sqrt{n \log n})$, but $\chi(G) \geq n/2$.}.

\section{\texorpdfstring{A constant-factor improvement for dominating models}{A constant-factor improvement for dominating models}}
\label{s:uniform_improvement}

In this section, we prove \cref{thm:uniform_improvement}, which is a constant-factor improvement of the Duchet--Meyniel bound for dominating models.

Let $S \subseteq V(G)$. We say that $S$ is \defn{connected} if $G[S]$ is connected. The set $S$ is \defn{dominating} if each vertex in $V(G) \setminus S$ has a neighbour in $S$. Write $\alpha(S):=\alpha(G[S])$, and let \defn{$c(S)$} denote the number of connected components of $G[S]$. Define the \defn{potential} of $S$ to be $\mathdefn{\phi(S)}:= 2\alpha(S) - c(S) - |S|$. The potential of a claw is one. A useful property of $\phi$ is that it is additive on disjoint vertex subsets with no edges between them. This follows from the fact that the independence number,
the number of components, and the number of vertices are all additive on disjoint vertex subsets with no edges between them. This notion is due to \citet{Fox10}.

A graph $G$ is \defn{$\alpha$-connected} if there is no partition of $V(G)$ into sets $A$ and $B$ with $A,B\neq\varnothing$ such that $\alpha(A)+\alpha(B)=\alpha(G)$.\footnote{This is the notion of `non-decomposability' in \citet{BK11} and \citet{Fox10}, but we use the term $\alpha$-connected since decomposable is already defined.}

\begin{thm}[{\citet[Lemma 2.1]{BK11}}]
    \label{thm:potential_cds}
    Let $G$ be an $n$-vertex $\alpha$-connected graph with $n\geqslant1$. For every $X \subseteq V(G)$, $G$ contains a connected dominating set $\tilde{X}$ with $|\tilde{X}| \leq 2\alpha(G) - 2\phi(X)/3 - 1$.
\end{thm}
We now prove \cref{thm:uniform_improvement}.
The argument does not extend to odd models, because its final step uses a connected dominating set as a new first branch set. We highlight this in the proof.

\uniformimprovementthm*
The proof follows the potential method of \citet{BK11} and \citet{Fox10}. We first choose a maximal collection of $r$ pairwise vertex-disjoint induced claws with no edges between distinct claws. Let $A$ be their union and $B:= N_G[A]$ its closed neighbourhood. After removing $B$, the remaining graph $D:= G - B$ is claw-free, so \cref{thm:odd_dom-claw-free} implies that $\haddom(D)$ is sufficiently large if $|V(D)|$ is large, and thus we are done in that case. This forces $|V(D)|$ to be small and hence $|B|$ to be large. If $\alpha(B)$ is small, then we apply induction on $B$ to deduce that $\haddom(G[B])$ is large. Thus $\alpha(B)$ is large and $B$ has a large independent set $I$. We then argue that $A \cup I$ is a set of large potential. Applying \cref{thm:potential_cds} to $A \cup I$, we obtain a small connected dominating set of $G$, and finish by induction, using it as the first branch set of the dominating model.
The main difference from \citet{BK11} is that their proof requires a more complicated analysis to find a set of large potential, involving several cases according to the size of $I$ and the structure of $D$. We avoid this analysis, at the cost of a slightly worse constant.

\begin{proof}
    Note that $0<c<1/9$ and $36c^2-49c+2=0$. Let $G$ be a minimal counterexample to \cref{thm:uniform_improvement}, and let $n:=|V(G)|$. Thus $n>(2-c)\alpha(G)\haddom(G)$.

    If $V(G)$ had a partition into sets $A$ and $B$ with $A,B\neq\varnothing$ and $\alpha(A)+\alpha(B)=\alpha(G)$, then by minimality,
    \begin{equation}
        |A|\leq (2-c)\alpha(A)\haddom(G[A])\leq (2-c)\alpha(A)\haddom(G). \label{eq:alpha-connected-part}
    \end{equation}
    Together with the analogous inequality for $B$, \cref{eq:alpha-connected-part} gives $n\leq (2-c)\alpha(G)\haddom(G)$, a contradiction. Hence $G$ is $\alpha$-connected. In particular, $G$ is connected.

    If $G$ is claw-free, then by \cref{thm:odd_dom-claw-free}, $\haddom(G)\geq 2\chi(G)/3\geq 2n/(3\alpha(G))>n/((2-c)\alpha(G))$, where the final inequality follows from $c<1/2$, a contradiction.

    Thus $G$ has an induced claw. Choose a maximal family of $r$ pairwise disjoint induced claws with no edges between distinct claws. Let $A$ be the union of their vertex sets, let $B:=N_G[A]$ be the closed neighbourhood, and let $D:=G-B$. Then
    \begin{equation*}
        |A|=4r,\qquad \alpha(A)=3r,\qquad c(A)=r,
        \qquad \phi(A)=r.
    \end{equation*}
    \textbf{Claim.}
    $G$ contains a set with potential at least $3c\alpha(G)/2$.
    \begin{proof}
        If $r\geq 3c\alpha(G)/2$, then $A$ is such a set. Thus assume that $r<3c\alpha(G)/2$.
        Let $I$ be a maximum independent set of $G[B]$, and consider $A\cup I$. This set has independence number $\alpha(B)$, at most $r$ connected components since each vertex of $I\setminus A$ is adjacent to a vertex of $A$, and at most $4r+\alpha(B)$ vertices. Hence
        \begin{equation*}
            \phi(A\cup I)\geq\alpha(B)-5r.
        \end{equation*}
        If $D$ is empty, then $\alpha(B)=\alpha(G)$, and since $c<1/9$ and $r<3c\alpha(G)/2$,
        \begin{equation*}
            \alpha(B)=\alpha(G)>9c\alpha(G)>5r+\frac{3c\alpha(G)}{2}.
        \end{equation*}
        Thus $A\cup I$ has the required potential. Hence assume that $D$ is non-empty, so $B\neq V(G)$.

        The maximality of $r$ implies that $D$ is claw-free. Since there are no edges between $A$ and $D$, we also have $\alpha(D)\leq\alpha(G)-3r$. By \cref{thm:odd_dom-claw-free},
        \begin{equation}
            |V(D)|\leq \frac{3}{2}\alpha(D)\haddom(D)
            \leq \frac{3}{2}(\alpha(G)-3r)\haddom(G). \label{eq:D-bound}
        \end{equation}
        Since $G$ is a counterexample, $\haddom(G)<n/((2-c)\alpha(G))$. Combining this with \cref{eq:D-bound} gives
        \begin{equation}
            |V(D)|\leq\frac{3(\alpha(G)-3r)n}{2(2-c)\alpha(G)}. \label{eq:D-size}
        \end{equation}
        Since $V(G)=B\cup V(D)$, \cref{eq:D-size} gives
        \begin{equation}
            |B|\geq n-\frac{3(\alpha(G)-3r)n}{2(2-c)\alpha(G)}
            =\frac{((1-2c)\alpha(G)+9r)n}{2(2-c)\alpha(G)}. \label{eq:B-size}
        \end{equation}
        Since $|B| < |V(G)|$, $G[B]$ satisfies \cref{thm:uniform_improvement}, so
        \begin{equation}
            |B|\leq (2-c)\alpha(B)\haddom(G)<\frac{\alpha(B)n}{\alpha(G)}. \label{eq:B-minimality}
        \end{equation}
        Combining \cref{eq:B-size,eq:B-minimality} gives
        \begin{equation}
            \alpha(B)>\frac{|B|\alpha(G)}{n}
            \geq\frac{(1-2c)\alpha(G)+9r}{2(2-c)}. \label{eq:alpha-B-lower}
        \end{equation}
        Since $36c^2-49c+2=0$, we have $1-8c+3c^2=3c(11-10c)/2$. Therefore
        \begin{equation}
            \begin{aligned}
                 & \frac{(1-2c)\alpha(G)+9r}{2(2-c)}-\left(5r+\frac{3c\alpha(G)}{2}\right)                                             \\
                 & \qquad=\frac{(1-8c+3c^2)\alpha(G)-(11-10c)r}{2(2-c)}=\frac{(11-10c)\left(\frac{3c\alpha(G)}{2}-r\right)}{2(2-c)}>0.
            \end{aligned}
            \label{eq:alpha-B-comparison}
        \end{equation}
        The final inequality follows from $r<3c\alpha(G)/2$ and $11-10c>0$. By \cref{eq:alpha-B-lower,eq:alpha-B-comparison}, we have $\alpha(B)>5r+3c\alpha(G)/2$. Thus $\phi(A\cup I)>3c\alpha(G)/2$, proving the claim.
    \end{proof}
    By \cref{thm:potential_cds}, applied to such a set, we deduce that $G$ has a connected dominating set $S$ with size at most
    \begin{equation}
        |S| \leq 2\alpha(G)-\frac{2}{3}\left(\frac{3c\alpha(G)}{2}\right)-1=(2-c)\alpha(G)-1 \label{eq:cds-size}.
    \end{equation}
    Since $n>(2-c)\alpha(G)\haddom(G)\geq (2-c)\alpha(G)$, $G-S$ has at least one vertex. Let $T$ be a spanning tree of $G[S]$. Prepending\footnote{This is the only place where the extension to odd models fails: even if $G[S]$ is connected and $S$ is dominating, we cannot conclude that $\hadodd(G)\geq \hadodd(G-S)+1$.} $T$ to a dominating model in $G - S$ shows that $\haddom(G - S) \leq \haddom(G) - 1$.
    By minimality,
    \begin{equation}
        |V(G - S)| \leq (2-c)\alpha(G)(\haddom(G)-1). \label{eq:remaining-size}
    \end{equation}
    By \cref{eq:cds-size,eq:remaining-size},
    \begin{equation*}
        n\leq (2-c)\alpha(G)-1+(2-c)\alpha(G)(\haddom(G)-1)
        =(2-c)\alpha(G)\haddom(G)-1,
    \end{equation*}
    a contradiction.
\end{proof}

\section{Odd dominating models}\label{sec:A}

The proof of \cref{thm:A} produces odd dominating models with the following additional property. A dominating $K_t$-model $(T_1,\dots,T_t)$ together with a colouring $c$ of $\bigcup_iV(T_i)$ with colours $\{\text{red},\text{blue}\}$ is \defn{red-dominating} if every edge of every $T_i$ is bichromatic, every $T_i$ contains a vertex of colour red, and for $1\leqslant i<j\leqslant t$ every vertex of $T_j$ has a neighbour of colour red in $T_i$. Every red-dominating $K_t$-model is an odd dominating $K_t$-model, since for $i<j$ a vertex of colour red in $T_j$ has a neighbour of colour red in $T_i$. Let \defn{$\hadred(G)$} be the maximum integer $t$ such that $G$ has a red-dominating $K_t$-model. Then $\hadred(G)\leqslant\hadod(G)$, so \cref{thm:A} is a consequence of the following result.

\begin{thm}\label{thm:Aplus}
    For every graph $G$ with $n\geqslant1$ vertices, there is an integer $h$ with $1\leqslant h\leqslant\hadred(G)$ such that $G$ has an $h$-colourable induced subgraph on at least $(n+h)/2$ vertices. In particular,
    \[ \hadred(G)\ \geqslant\ \frac{n}{2\alpha(G)-1}. \]
\end{thm}
The following lemma is implicit in the proofs of \citet{KS07,DM82}.

\begin{lem}\label{lem:A}
    Every non-empty connected graph $G$ has a set $D\subseteq V(G)$, an independent set $I\subseteq D$ with $|D|\leqslant 2|I|-1$, and a spanning tree $F$ of $G[D]$ such that
    \begin{enumerate}[label=(\arabic*)]
        \item every edge of $F$ has exactly one end in $I$, and
        \item every vertex of $V(G)\setminus D$ has a neighbour in $I$.
    \end{enumerate}
\end{lem}

\begin{proof}
    There is at least one triple $(D, I, F)$ satisfying all the stated conditions except possibly ~(2).
    For example, take $D=I=\{v\}$ for any vertex $v$, and let $F$ be the 1-vertex tree. Choose such a triple with $|I|$ maximum, and suppose that (2) fails. Let $N$ be the set of vertices in $V(G)\setminus D$ with a neighbour in $I$, and let $R:=V(G)\setminus(D\cup N)$. Then $R\neq\varnothing$. Since $G$ is connected, there is an edge $yz$ with $y\in R$ and $z\notin R$. Since $y\notin N$, the vertex $y$ has no neighbour in $I$, so $z\notin I$. Thus $z\in N$ or $z\in D\setminus I$.
    \begin{itemize}
        \item If $z\in N$, then let $x\in I$ be a neighbour of $z$, and let $D':=D\cup\{y,z\}$, $I':=I\cup\{y\}$ and $F':=F+xz+zy$. Then $I'$ is independent, each edge of $F'$ has exactly one end in $I'$, and $|D'|=|D|+2\leqslant 2|I'|-1$.
        \item If $z\in D\setminus I$, then let $D':=D\cup\{y\}$, $I':=I\cup\{y\}$ and $F':=F+zy$. Then $I'$ is independent, each edge of $F'$ has exactly one end in $I'$, and $|D'|=|D|+1<2|I'|-1$.
    \end{itemize}
    In both cases $|I'|>|I|$, which contradicts the choice of $(D,I,F)$.
\end{proof}

In \cref{lem:A}, colour each vertex in $I$ red and each vertex in $D\setminus I$ blue. Then every edge of $F$ is bichromatic, and every vertex of $V(G)\setminus D$ has a neighbour of colour red in $D$. The analogous lemma of \citet{IW24} only guarantees that every vertex of $V(G)\setminus D$ has a neighbour in $D$. This is not enough for oddness, since all the neighbours in $D$ of a vertex outside $D$ might have colour blue.

For a rooted forest $T$ and vertices $x,y\in V(T)$, we say that $x$ is a \defn{proper ancestor} of $y$ if $x\neq y$ and $x$ lies on the path in $T$ from $y$ to the root of its component. Vertices $x,y$ are \defn{comparable} if $x=y$ or one is a proper ancestor of the other. Applying \cref{lem:A} recursively gives the following structure theorem.

\begin{thm}\label{thm:struct}
    For every graph $G$ there is a rooted forest $T$ and a partition $\{P_x:x\in V(T)\}$  of $V(G)$ into non-empty parts, together with, for each $x\in V(T)$, an independent set $I_x\subseteq P_x$ with $|P_x|\leqslant 2|I_x|-1$ and a spanning tree $F_x$ of $G[P_x]$ whose edges each have exactly one end in $I_x$, such that:
    \begin{enumerate}[label=(\arabic*)]
        \item if $x,y\in V(T)$ are incomparable, then there is no edge of $G$ between $P_x$ and $P_y$;
        \item if $x$ is a proper ancestor of $y$, then every vertex of $P_y$ has a neighbour in $I_x$.
    \end{enumerate}
\end{thm}

\begin{proof}
    We proceed by induction on $|V(G)|$. If $G$ has no vertices, let $T$ be empty. If $G$ is disconnected, then apply induction to each component of $G$ and take the union of the resulting forests. Then (1) holds since distinct components are not adjacent, and (2) holds within each component. Now assume that $G$ is connected. Apply \cref{lem:A} to obtain $(D,I,F)$. By induction, $G-D$ has a rooted forest $T'$ and a partition $\{P_x:x\in V(T')\}$ with the stated properties. Let $T$ be obtained from $T'$ by adding a new vertex $r$ adjacent to the root of each component of $T'$, and let $r$ be the root of $T$. Let $P_r:=D$, $I_r:=I$ and $F_r:=F$. Property~(1) holds since $r$ is comparable with every vertex of $T$, and incomparable vertices of $T$ are incomparable in $T'$. Property~(2) holds for $x=r$ by \cref{lem:A}(2), and for $x\neq r$ by induction.
\end{proof}

Paths from a root in this structure give red-dominating models.

\begin{cor}\label{cor:chain}
    In the setting of \cref{thm:struct}, let $c(v):=\text{red}$ if $v\in\bigcup_xI_x$ and $c(v):=\text{blue}$ otherwise. If $x_1$ is a root of $T$ and $x_{i+1}$ is a child of $x_i$ for $1\leqslant i<h$, then $((F_{x_1},\dots,F_{x_h}),c)$ is a red-dominating $K_h$-model in $G$.
\end{cor}

\begin{proof}
    Each $F_{x_i}$ is a tree whose edges are bichromatic, and $F_{x_i}$ contains the vertices of $I_{x_i}\neq\varnothing$, which have colour red. For $i<j$, the vertex $x_i$ is a proper ancestor of $x_j$, so every vertex of $F_{x_j}$ has a neighbour in $I_{x_i}$, which has colour red.
\end{proof}

\begin{proof}[Proof of \cref{thm:Aplus}]
    Apply \cref{thm:struct} to $G$, and let $h$ be the maximum number of vertices on a path in $T$ that starts at a root. Since $G$ has at least one vertex, $h\geqslant1$, and $h\leqslant\hadred(G)$ by \cref{cor:chain}. Let $X:=\bigcup_{x\in V(T)}I_x$. Then
    \[ |X|=\sum_{x\in V(T)}|I_x|\ \geqslant\ \sum_{x\in V(T)}\frac{|P_x|+1}{2}=\frac{n+|V(T)|}{2}\ \geqslant\ \frac{n+h}{2}. \]
    Colour each vertex of $I_x$ by the depth of $x$ in $T$, where roots have depth~1. This uses at most $h$ colours. Two distinct vertices of $X$ with the same colour are either both in the same independent set $I_x$, or in parts $P_x$ and $P_y$ where $x\neq y$ have the same depth. In the latter case, $x$ and $y$ are incomparable, so the two vertices are non-adjacent by \cref{thm:struct}(1). Thus $G[X]$ is $h$-colourable. Some colour class of $G[X]$ is an independent set of size at least
    \[ \frac{n+h}{2h}=\frac{n}{2h}+\frac12\ \geqslant\ \frac{n}{2\hadred(G)}+\frac12=\frac{n+\hadred(G)}{2\hadred(G)}. \]
    Thus $\alpha(G)\geqslant(n+\hadred(G))/(2\hadred(G))$, which rearranges to $\hadred(G)\geqslant n/(2\alpha(G)-1)$.
\end{proof}

The next result bounds $\chi(G)$ in terms of $\hadred(G)$. The proof adapts Wagner's argument~\citep{Wag64} for graphs with no $K_t$-model.

\begin{prop}\label{prop:wagner}
    Every graph $G$ with at least one vertex satisfies $\chi(G)\leqslant 2^{\hadred(G)-1}$. In particular, $\chi(G)\leqslant2^{\hadod(G)-1}$.
\end{prop}

\begin{proof}
    We proceed by induction on $k \geq 1$ to show that every non-empty graph $G$ with $\hadred(G)\leqslant k$ satisfies
    $\chi(G)\leqslant2^{k-1}$. If $k=1$, then $|E(G)| = 0$, so $\chi(G)=1$.

    Now let $k\geqslant 2$. We may assume that $G$ is connected.
    Fix a vertex $r$, let $L_0:= \{r\}$, and for $i \geq 1$, let $L_i$ be the set of vertices at distance exactly $i$ from $r$.
    Every edge of $G$ has both ends in one $L_i$ or ends in $L_i$ and $L_{i+1}$ for some $i$.
    Consequently the sets $L_i$ with $i$ even are pairwise non-adjacent, as are the sets $L_i$ with $i$ odd, and thus $\chi(G)\leqslant2\max_i\chi(G[L_i])$. We now show that $\chi(G[L_i]) \leq 2^{k-2}$ for every layer $L_i$, which implies the result.

    This holds for $L_0$ since $k\geqslant2$.
    For a layer $L_i$ with $i\geqslant1$, take a red-dominating $K_s$-model in $G[L_i]$, with  $s:=\hadred(G[L_i])$. Let $T$ be a breadth-first search tree from $r$ restricted to $L_0\cup\dots\cup L_{i-1}$, which is a subtree of $G$. Colour each vertex of $T$ at distance $d$ from $r$ with red if $d\equiv i-1\pmod 2$ and blue otherwise. Every edge of $T$ is bichromatic, and every vertex of $L_{i-1}$ is red. Since every vertex of $L_i$ has a
    neighbour in $L_{i-1}$, prepending $T$ to the red-dominating $K_s$ model gives
    a red-dominating $K_{s+1}$-model in $G$. Thus $s+1 \leq \hadred(G) \leq k$. This implies that $s = \hadred(G[L_i]) \leq k - 1$. By induction, $\chi(G[L_i]) \leq 2^{k-2}$, completing the proof.
\end{proof}

\section{Strongly odd dominating models}\label{sec:B}

In a strongly odd dominating model, every vertex of a later branch set must have a neighbour of its own colour in each earlier branch set.
The following strengthening of \cref{lem:A} ensures the stronger property that every vertex outside the constructed set has neighbours of both colours in it.

\begin{lem}\label{lem:B}
    Every non-empty connected graph $G$ has a set $D\subseteq V(G)$, a colouring $c:D\to\{\text{red},\text{blue}\}$, a spanning tree $F$ of $G[D]$ whose edges are bichromatic, and an independent set $Y\subseteq D$ with $|Y|\geqslant(|D|+1)/3$, such that every vertex of $V(G)\setminus D$ has a neighbour of colour red and a neighbour of colour blue in $D$.
\end{lem}

\begin{proof}
    Let $(D_0,I,F_0)$ be given by \cref{lem:A}. Colour $I$ with colour red and $D_0\setminus I$ with colour blue. Let $W$ be the set of vertices in $V(G)\setminus D_0$ with no neighbour in $D_0\setminus I$, and let $J$ be a maximal independent set of $G[W]$, where $J=\varnothing$ if $W=\varnothing$. By \cref{lem:A}(2), every $w\in J$ has a neighbour $x_w\in I$. Let $D:=D_0\cup J$, colour $J$ with colour blue, and let $F:=F_0+\{x_ww:w\in J\}$. Then $F$ is a spanning tree of $G[D]$ whose edges are bichromatic.

    Let $v\in V(G)\setminus D$. By \cref{lem:A}(2), $v$ has a neighbour in $I$, which has colour red. If $v\notin W$, then $v$ has a neighbour in $D_0\setminus I$, which has colour blue. If $v\in W\setminus J$, then $v$ has a neighbour in $J$ by the maximality of $J$, and this neighbour has colour blue.

    Finally, $|D|=|D_0|+|J|\leqslant 2|I|+|J|-1\leqslant 3\max\{|I|,|J|\}-1$. Let $Y$ be the larger of $I$ and $J$. Then $Y$ is independent and $|Y|\geqslant(|D|+1)/3$.
\end{proof}

Using \cref{lem:B} in place of \cref{lem:A} gives the following analogue of \cref{thm:struct}.
\begin{thm}\label{thm:structB}
    For every graph $G$ there is a rooted forest $T$, a partition $\{P_x:x\in V(T)\}$ of $V(G)$ into non-empty parts, and a colouring $c:V(G)\to\{\text{red},\text{blue}\}$, together with, for each $x\in V(T)$, a spanning tree $F_x$ of $G[P_x]$ whose edges are bichromatic and an independent set $Y_x\subseteq P_x$ with $|Y_x|\geqslant(|P_x|+1)/3$, such that:
    \begin{enumerate}[label=(\arabic*)]
        \item if $x,y\in V(T)$ are incomparable, then there is no edge of $G$ between $P_x$ and $P_y$;
        \item if $x$ is a proper ancestor of $y$, then every vertex of $P_y$ has a neighbour of colour red and a neighbour of colour blue in $P_x$.
    \end{enumerate}
\end{thm}

\begin{proof}
    The proof is the same as that of \cref{thm:struct}, with \cref{lem:B} in place of \cref{lem:A}. For connected $G$, the root part is the set $D$ given by \cref{lem:B}, with its colouring, tree $F$ and independent set $Y$, and the colouring of $V(G)\setminus D$ is given by induction.
\end{proof}

We now prove \cref{thm:B}, following the proof of \cref{thm:Aplus} with \cref{thm:structB} in place of \cref{thm:struct}.

\begin{proof}[Proof of \cref{thm:B}]
    Apply \cref{thm:structB} to $G$, and let $h$ be the maximum number of vertices on a path in $T$ that starts at a root. Since $G$ has at least one vertex, $h\geqslant1$. If $x_1,\dots,x_h$ is such a path with $x_1$ a root, then $((F_{x_1},\dots,F_{x_h}),c)$ is a strongly odd dominating $K_h$-model, since for $i<j$ each vertex $v\in P_{x_j}$ has a neighbour of colour $c(v)$ in $P_{x_i}$ by \cref{thm:structB}(2). Thus $h\leqslant\hadsod(G)$. Let $X:=\bigcup_{x\in V(T)}Y_x$. Then
    \[ |X|=\sum_{x\in V(T)}|Y_x|\ \geqslant\ \sum_{x\in V(T)}\frac{|P_x|+1}{3}=\frac{n+|V(T)|}{3}\ \geqslant\ \frac{n+h}{3}. \]
    As in the proof of \cref{thm:Aplus}, colouring each vertex of $Y_x$ by the depth of $x$ in $T$ shows that $G[X]$ is $h$-colourable. Some colour class is an independent set of size at least
    \[ \frac{n+h}{3h}=\frac{n}{3h}+\frac13\ \geqslant\ \frac{n}{3\hadsod(G)}+\frac13=\frac{n+\hadsod(G)}{3\hadsod(G)}. \]
    Thus $\alpha(G) \geqslant (n+\hadsod(G))/(3\hadsod(G))$, which rearranges to $\hadsod(G) \geqslant n/(3\alpha(G)-1)$.
\end{proof}

We now bound $\chi(G)$ in terms of $\hadsod(G)$. The proof of \cref{prop:wagner} does not work here, since in a breadth-first search all the neighbours of a vertex of $L_i$ in $L_0\cup\dots\cup L_{i-1}$ lie in $L_{i-1}$, and so receive the same colour in that proof. The following lemma deals with this issue.

\begin{lem}\label{lem:absorb}
    Let $W\subseteq V(G)$ be a non-empty set such that $G[W]$ is connected, let $F$ be a spanning tree of $G[W]$, let $c:W\to\{\text{red},\text{blue}\}$ be a colouring such that every edge of $F$ is bichromatic, and let $K\subseteq V(G)\setminus W$ be a set such that every vertex of $K$ has a neighbour in $W$. Then there is a set $Z\subseteq K$, an extension of $c$ to $W\cup Z$, and a spanning tree $F'$ of $G[W\cup Z]$ containing $F$, such that:
    \begin{enumerate}[label=(\arabic*)]
        \item every edge of $F'$ is bichromatic;
        \item every edge of $G[Z]$ is bichromatic, so $\chi(G[Z])\leqslant2$;
        \item every vertex of $K\setminus Z$ has a neighbour of colour red and a neighbour of colour blue in $W\cup Z$.
    \end{enumerate}
\end{lem}

\begin{proof}
    Start with $Y:=W$, $Z:=\varnothing$ and $F':=F$. While some vertex $v\in K\setminus Z$ has neighbours in $Y$ of only one colour $\gamma$ (it has at least one neighbour in $Y$, since it has one in $W$), add $v$ to $Y$ and to $Z$, colour $v$ with the colour in $\{\text{red},\text{blue}\}\setminus\{\gamma\}$, and add to $F'$ an edge from $v$ to a neighbour of $v$ of colour $\gamma$. This process terminates, and then (1) and (3) hold. For (2), let $zz'$ be an edge of $G[Z]$, where $z$ was added before $z'$. When $z'$ was added, $z$ was in $Y$ and adjacent to $z'$, so the only colour that $z'$ saw in $Y$ was the colour of $z$, and $z'$ received the other colour.
\end{proof}

With \cref{lem:absorb}, the layering argument gives the following analogue of \cref{prop:wagner}.

\begin{prop}\label{prop:sod}
    Every graph $G$ with at least one vertex satisfies $\chi(G)\leqslant3\cdot2^{\hadsod(G)}-5$.
\end{prop}

\begin{proof}
    We proceed by induction on $k\geqslant1$ to show that every
    non-empty graph $G$ with $\hadsod(G)\leqslant k$ satisfies
    $\chi(G)\leqslant3\cdot2^k-5$. If $k=1$, then $|E(G)|=0$,
    so $\chi(G)=1$.

    Now let $k\geqslant2$. We may assume that $G$ is connected.
    Fix a vertex $r$, let $L_0:=\{r\}$, and for $i\geqslant1$,
    let $L_i$ be the set of vertices at distance exactly $i$
    from $r$.
    Every edge of $G$ has both ends in one $L_i$ or ends in $L_i$ and $L_{i+1}$ for some $i$.
    Consequently the sets $L_i$ with $i$ even are pairwise non-adjacent, as are the sets $L_i$ with $i$ odd, and thus $\chi(G)\leqslant2\max_i\chi(G[L_i])$.
    We now show that $\chi(G[L_i])\leqslant3\cdot2^{k-1}-3$
    for every layer $L_i$, which implies the result.
    This holds for $L_0$ since $k\geqslant2$.

    For a layer $L_i$ with $i\geqslant1$,
    let $W:=L_0\cup\dots\cup L_{i-1}$, and let $F$ be a
    breadth-first search spanning tree from $r$ restricted to $W$.
    Let $c$ be the colouring of $W$ in which each vertex at distance
    $d$ from $r$ is red if $d\equiv i-1\pmod2$ and blue otherwise.
    Each edge of $F$ is bichromatic, and each vertex of $L_i$
    has a neighbour in $L_{i-1}\subseteq W$.
    Apply \cref{lem:absorb} to $W$, $F$, $c$ and $L_i$
    to obtain $Z$, the extension of $c$ (also denoted $c$)
    and a spanning tree $F'$. Let $H:=G[L_i\setminus Z]$.
    By \cref{lem:absorb}(2),
    $\chi(G[L_i])\leqslant\chi(H)+2$.
    If $H$ is empty, then $\chi(G[L_i])\leqslant2$, as required.
    Otherwise, take a strongly odd dominating $K_s$-model
    $((T_1,\dots,T_s),c')$ in $H$, with $s:=\hadsod(H)$.
    By \cref{lem:absorb}(3), for each $j\in\{1,\dots,s\}$,
    every vertex of $T_j$ has neighbours of both colours in $F'$,
    and in particular a neighbour of its own colour under $c'$.
    Hence $((F',T_1,\dots,T_s),c\cup c')$ is a strongly odd
    dominating $K_{s+1}$-model in $G$. Thus $s+1\leqslant\hadsod(G)\leqslant k$, so
    $\hadsod(H)\leqslant k-1$. By induction,
    $\chi(H)\leqslant3\cdot2^{k-1}-5$. Consequently, $\chi(G[L_i])\leqslant3\cdot2^{k-1}-3$,
    completing the proof.
\end{proof}

\section{Open Problems}

We finish with some open problems.
\begin{enumerate}[label=(\arabic*)]
    \item The best known upper bound is $\chi(G)\in O(\haddom(G)(\log\haddom(G))^2)$ for every graph $G$ with $\haddom(G)\geqslant2$, by the result of \citet{GNTT26} mentioned above. Is there a constant $C$ such that every graph $G$ satisfies $\chi(G)\leqslant C\haddom(G)$? \citet{IW24} showed that there are graphs with no dominating $K_t$-model and with average degree $(1-o(1))t\log t$, so a bound on the degeneracy does not give a linear bound on the chromatic number.
    \item Every non-empty graph $G$ satisfies $\chi(G)\leqslant2^{\hadod(G)-1}$ (\cref{prop:wagner}, which adapts Wagner's layering argument~\citep{Wag64}) and $\chi(G)\leqslant3\cdot2^{\hadsod(G)}-5$ (\cref{prop:sod}). Is $\chi(G)$ bounded above by a polynomial function of $\hadod(G)$? The same question applies to $\hadsod(G)$. Bipartite graphs $G$ satisfy $\hadodd(G)\leqslant2$ and can have arbitrarily large minimum degree, so the degeneracy of $G$ is not bounded by any function of $\hadod(G)$ or $\hadsod(G)$. Thus such a result cannot be proved by bounding the degeneracy.
    \item What are the constants $\rho$ and $\rho_d$?
\end{enumerate}

\paragraph{Use of AI.} The authors discovered all the proofs in this paper. Claude Opus 5.5 was used to help write the paper, and GPT-6.1 was used for proof-reading.

{
\fontsize{10pt}{11pt}
\selectfont
\bibliographystyle{DavidNatbibStyle}
\bibliography{ref}
}

\section*{\texorpdfstring{Appendix: Proofs of \cref{thm:dominating_seagulls,thm:plummer}}{Appendix: Seagull and clique-number proofs}}
\label{sec:appendix}

The proof of \cref{thm:dominating_seagulls} extends a deep result of \citet[1.3]{CS12} to odd dominating models.

A \defn{seagull} is an induced three-vertex path. A \defn{five-wheel} consists of a $5$-cycle and one additional vertex adjacent to every vertex of the cycle.
For disjoint sets $X,Y\subseteq V(G)$, we say that $X$ is
\emph{complete to} $Y$ if every vertex of $X$ is adjacent to every
vertex of $Y$, and \emph{anticomplete to} $Y$ if there are no edges
between $X$ and $Y$.
For a non-empty
clique $C$ of a graph $G$, let $A_C$, $B_C$, and $D_C$ denote the
sets of vertices in $V(G)\setminus C$ that are complete to $C$,
anticomplete to $C$, and have both a neighbour and a non-neighbour
in $C$, respectively. We call $(A_C, B_C, C, D_C)$ the \defn{associated partition} of a clique $C$. Define the \defn{capacity} of $C$ by
\[
    \operatorname{cap}_G(C)
    := |D_C|+\frac{|A_C|+|B_C|}{2}.
\]

To prove \cref{thm:dominating_seagulls}, the following result of \citet{CS12} is essential.

\begin{thm}[{\citet[1.6]{CS12}}]
    \label{thm:seagulls}
    Let $G$ be a non-empty graph with $\alpha(G)\leq 2$, and let $k\geq 0$ be an
    integer, such that if $k=2$ then $G$ is not a five-wheel.
    Then $G$ has $k$ pairwise vertex-disjoint seagulls if and only if
    \begin{itemize}
        \item $|V(G)|\geq 3k$;
        \item $G$ is $k$-connected;
        \item every non-empty clique of $G$ has capacity at least $k$; and
        \item $\overline{G}$ has a matching of cardinality $k$.
    \end{itemize}
\end{thm}
We now prove \cref{thm:dominating_seagulls} using \cref{thm:seagulls}.

\dominatingseagulls*

\begin{proof}
    Let $C$ be a clique of maximum size in $G$, and we may suppose that $|C|\leq \ceil{n/2} - 1$. Let $H:= G - C$. Let $k:=\lceil n/2\rceil-|C| \geq 1$.

    \textbf{Claim.} $H$ satisfies the hypotheses of \cref{thm:seagulls} and thus has $k$ pairwise vertex-disjoint seagulls.

    First, we check that $|V(H)|\geq 3k$.
    By hypothesis,
    $2|C|\geq \ceil{n/2}+1$, and $2\ceil{n/2}\leq n+1$. Hence
    \begin{align*}
        |V(H)|
         & =n-|C|                                         \\
         & =3\bigl(\ceil{n/2}-|C|\bigr)
        +n+2|C|-3\ceil{n/2}                               \\
         & =3k+n+2|C|-3\ceil{n/2}                         \\
         & \geq 3k+n+\bigl(\ceil{n/2}+1\bigr)-3\ceil{n/2} \\
         & =3k+n+1-2\ceil{n/2}\geq 3k.
    \end{align*}
    Next, we show that $H$ is $k$-connected. Suppose for contradiction that there is a set
    $X\subseteq V(H)$ with $|X|\leq k-1$ such that $H-X$ is
    disconnected. As $\alpha(G)\leq 2$, the graph $H-X$ consists of exactly two components, each of which is a clique.
    Let their vertex sets be $A$ and $B$.

    Now, each vertex  $v\in C$ is complete to either $A$ or $B$ (otherwise taking a vertex $a \in A$ and $b \in B$ not adjacent to $v$ gives an independent set of size $3$). Thus we can partition $C$ into sets
    $C_A$ and $C_B$ such that $C_A$ is complete to $A$ and
    $C_B$ is complete to $B$. Both $A\cup C_A$ and $B\cup C_B$
    are cliques, so the maximality of $|C|$ implies that
    \begin{equation*}
        n-|X|
        =|A\cup C_A|+|B\cup C_B|
        \leq 2|C|.
    \end{equation*}
    On the other hand, since $n\geq 2\ceil{n/2}-1$,
    \begin{equation*}
        n-|X|
        \geq n-k+1
        =n-\ceil{n/2}+|C|+1
        \geq \ceil{n/2}+|C|
        >2|C|,
    \end{equation*}
    a contradiction.

    Next, we show that every non-empty clique $Q$ of $H$ satisfies $\operatorname{cap}_H(Q)\geq k$.
    Suppose that a non-empty clique $Q$ of $H$ satisfies
    $\operatorname{cap}_H(Q)<k$.
    Let $(Q_A,Q_B,Q,Q_D)$ be its associated
    partition in $H$.
    Then since $|Q_A|, |Q_B|$ and $|Q_D|$ are integers,
    \begin{equation*}
        \operatorname{cap}_H(Q)
        =|Q_D|+\frac{|Q_A|+|Q_B|}{2} \leq k-\frac{1}{2}.
    \end{equation*}
    Therefore,
    \begin{equation*}
        |Q_A| + |Q_B| + 2|Q_D| \leq 2k - 1.
    \end{equation*}
    Since $2\ceil{n/2} - 1 \leq n$, we have
    \begin{align*}
        2\ceil{n/2} - 1 + |Q_D| & \leq n + |Q_D|                                       \\
                                & = |Q_A| + |Q_B| + |Q| + 2|Q_D| + |C|                 \\
                                & \leq 2k - 1 + |Q| + |C|                              \\
                                & \leq 2k - 1 + 2|C| = 2\left\lceil n/2\right\rceil-1,
    \end{align*}
    where the last inequality follows from $k=\left\lceil n/2\right\rceil-|C|$. Hence every inequality is an equality. In particular, $Q_D = \varnothing$ and $|Q| = |C|$. Thus $Q$ is a maximum clique of $G$ and $Q_A = \varnothing$. Thus $V(H) = Q \cup Q_B$. Since $H$ is $k$-connected and $k \geq 1$, $H$ is connected. In particular, $Q_B = \varnothing$. Consequently, $V(H) = Q$ and $V(G) = C \cup Q$. But $|Q| = |C| < n/2$, a contradiction. Thus every non-empty clique $Q$ of $H$ satisfies $\operatorname{cap}_H(Q)\geq k$.

    Next, we show that $\overline{H}$ has a matching of cardinality $k$. Let $M$ be a maximal matching in $\overline{H}$, and let $U$ be the set of vertices of $H$ not covered by $M$. By the maximality of $M$, the set $U$ is independent in $\overline{H}$, and hence is a clique in $H$. Thus, $|U| \leq |C|$. Therefore,
    \begin{align*}
        2|M|
         & = |V(H)|-|U|                                     \\
         & \geq n-|C|-|C|                                   \\
         & = n-2|C|                                         \\
         & \geq 2\left\lceil n/2\right\rceil-1-2|C| = 2k-1.
    \end{align*}
    Since $2|M|$ is even, we have $2|M| \geq 2k$ and hence $|M| \geq k$, as desired.

    Finally, to apply \cref{thm:seagulls}, we need to show that $H$ is not isomorphic to the $5$-wheel when $k = 2$. Suppose that $H$ is isomorphic to the $5$-wheel. Then $|V(H)| = 6$ and $H$ contains a clique of size $3$. Since $C$ is a maximum clique of $G$, we have $|C| \geq 3$. Thus $n \geq 9$. Since $k=2$, $n=|C|+6
        =\left\lceil n/2\right\rceil-k+6
        =\left\lceil n/2\right\rceil+4$. We deduce that $\floor{n/2} = 4$, and so $n \in \{8,9\}$. Consequently, $n = 9$ and $|C| = 3$.

    Let $D$ be the induced $5$-cycle in $H$. Since $C$ is a maximum clique of $G$, $\omega(G) \leq 3$. Therefore, for every edge $uv$ of $D$, at most one vertex of $C$ is adjacent to both $u$ and $v$. Since $|E(D)| = 5$ and $|C| = 3$, there is a vertex $c \in C$ such that $c$ is adjacent to both $u$ and $v$ for at most one edge $uv \in E(D)$.
    Thus $c$ has at most three neighbours in $D$, and the three neighbours are not consecutive on $D$. Let $x$ and $y$ be non-adjacent vertices in $D$ such that $x$ and $y$ are non-adjacent to $c$.
    But $\{x,y,c\}$ is an independent set of size $3$, a contradiction. Therefore, $H$ is not a $5$-wheel. By \cref{thm:seagulls}, $H$ contains $k$ pairwise vertex-disjoint seagulls, completing the proof of the claim.

    Let $S_1, \dots, S_k$ be the $k$ vertex-disjoint seagulls of $H$ supplied by the claim, and for each $i \in \{1, \dots, k\}$, let $S_i = x_i y_i z_i$, where $x_iz_i \notin E(H)$. Enumerate the vertices of $C = \{c_1, \dots, c_{|C|}\}$, and consider the $K_{k + |C|} = K_{\ceil{n/2}}$-model of $G$ given by the branch sets
    \begin{equation*}
        \mathcal{T}:=(S_1,\ldots,S_k,\{c_1\},\ldots,\{c_{|C|}\}).
    \end{equation*}
    We claim that $\mathcal{T}$ is an odd and dominating model. Since $\alpha(G)\leq 2$ and $x_i z_i\notin E(G)$, every vertex $v\in V(G)\setminus V(S_i)$ is adjacent to at least one of $x_i$ and $z_i$; otherwise, $\{v,x_i,z_i\}$ is an independent set of size $3$. Thus $\mathcal{T}$ is a dominating model.

    To see that $\mathcal{T}$ is an odd model, for each seagull $S_i$, we colour $x_i$ and $z_i$ blue and $y_i$ red. We colour each vertex of $C$ blue.
    The endpoints of every edge within each branch set receive different colours. Every blue vertex outside $S_i$ has a blue neighbour in $\{x_i,z_i\}$. Since $C$ is a blue clique, every two branch sets are joined by a monochromatic edge. Hence $\mathcal{T}$ is an odd model. This completes the proof.
\end{proof}

As mentioned, the proof of \cref{thm:plummer} is a straightforward modification of the proof in \citet[Theorem 3.3(a)]{PST03}.
\plummer*
\begin{proof}
    Let $C$ be a maximum clique of $G$, and let $S_1,\ldots,S_s$ be a maximal collection of $s \geq 0$ pairwise vertex-disjoint seagulls in $G-C$. Let $G':=G-\bigcup_{i=1}^s V(S_i)$ and $H:=G'-C$.
    By the maximality of the collection, $H$ contains no seagull. Thus every component of $H$ is a clique, and since $\alpha(G)\leq 2$, $H$ has at most two components; let $A$ and $B$ be their vertex sets, where $A$ or $B$ may be empty.

    Every vertex of $C$ is complete to $A$ or $B$, since otherwise there is an independent set of size $3$.
    Hence we may partition $C=C_A \cup C_B$ such that $C_A$ is complete to $A$ and $C_B$ is complete to $B$. Since $A\cup C_A$ and $B\cup C_B$ are cliques, $|V(G')|=|A\cup C_A|+|B\cup C_B|\leq 2|C|.$

    Since each seagull has three vertices, we have $n-3s=|V(G')|\leq 2|C|.$
    Therefore, $n+|C|\leq 3s+3|C|=3(s+|C|)$ and consequently $s+|C|\geq \ceil{(n + |C|)/3} = \ceil{(n + \omega(G))/3}$.

    For each $i\in\{1,\ldots,s\}$, write $S_i=x_i y_i z_i$, where $x_i z_i\notin E(G)$, and enumerate $C=\{c_1,\ldots,c_{|C|}\}$. Consider the model
    \begin{equation*}
        \mathcal{T}:=(S_1,\ldots,S_s,\{c_1\},\ldots,\{c_{|C|}\}).
    \end{equation*}
    We claim that $\mathcal{T}$ is an odd and dominating model. Since $\alpha(G)\leq 2$ and $x_i z_i\notin E(G)$, every vertex $v\in V(G)\setminus V(S_i)$ is adjacent to at least one of $x_i$ and $z_i$; otherwise, $\{v,x_i,z_i\}$ is an independent set of size $3$. Thus $\mathcal{T}$ is a dominating model.

    To see that $\mathcal{T}$ is an odd model, for each seagull $S_i$, we colour $x_i$ and $z_i$ blue and $y_i$ red. We colour each vertex of $C$ blue.
    The endpoints of every edge within each branch set receive different colours. Every blue vertex outside $S_i$ has a blue neighbour in $\{x_i,z_i\}$. Since $C$ is a blue clique, every two branch sets are joined by a monochromatic edge. Hence $\mathcal{T}$ is an odd model.
    Consequently, $\hadod(G)\geq s+|C|\geq \ceil{(n+\omega(G))/3},$ as required.
\end{proof}
\end{document}